\documentclass[12pt]{article}

\usepackage[margin=1in,nohead,a4paper,dvips]{geometry}
\usepackage{graphics}
\usepackage[center]{caption}

\makeatletter
\def\eqnarray{%
   \stepcounter{equation}%
   \def\@currentlabel{\p@equation\theequation}%
   \global\@eqnswtrue
   \m@th
   \global\@eqcnt\z@
   \tabskip\@centering
   \let\\\@eqncr
   $$\everycr{}\halign to\displaywidth\bgroup
       \hskip\@centering$\displaystyle\tabskip\z@skip{##}$\@eqnsel
      &\global\@eqcnt\@ne\hfil$\displaystyle{\hbox{}##\hbox{}}$\hfil
      &\global\@eqcnt\tw@ $\displaystyle{##}$\hfil\tabskip\@centering
      &\global\@eqcnt\thr@@ \hb@xt@\z@\bgroup\hss##\egroup
         \tabskip\z@skip
      \cr
}
\makeatother

\def\pairing #1#2{\scalebox{#2}{\includegraphics{p#1.eps}}}

\newtheorem{con}{Conjecture}

\begin{document}

\title{Combinatorial nature of ground state vector\\
of O(1) loop model}
\author{A.~V.~Razumov, Yu.~G.~Stroganov\\
\small \it Institute for High Energy Physics\\[-.5em]
\small \it 142280 Protvino, Moscow region, Russia}

\date{}

\maketitle

\begin{abstract}
Hanging about a hypothetical connections between the ground state vector
for some special spin systems and the alternating-sign matrices,
we have found a numerical evidence for the fact that the numbers of the
states of the fully packed loop model with fixed link-patterns
coincide with the components of the ground state vector of the dense O$(1)$
loop model considered by Batchelor, de Gier and Nienhuis. Our conjecture
generalizes in a sense the conjecture of Bosley and Fidkowski, refined by
Cohn and Propp, and proved by Wieland.

\end{abstract}

\leftskip 15em
{\it I  would  say  that imagination is a form of memory.}

\hbox to \hsize{\hfil (V.~Nabokov)}

\leftskip 0em

\vskip 1em

In paper~\cite{RS} we made some conjectures related to combinatorial
properties of the ground state vector of the XXZ spin chain for the
asymmetry parameter $\Delta=-1/2$ and an odd number of cites. In the
subsequent paper~\cite{BGN} Batchelor, de Gier and Nienhuis considered two
variations of this model along with the corresponding O$(n)$ loop model at
$n = 1$ and notably increased the number of models and related
combinatorial objects (see also \cite{RST}).

During last months we have accumulated a lot of data on relation
between the spin systems and alternating sign matrices (ASMs).
In this note we limit ourselves to one qualitative conjecture which relates
some classes of the states of the fully packed loop (FPL) model
with the ground state vector of the dense O$(n)$ loop model (see
\cite{BNi89} and references therein). The states the
FPL model is in bijective correspondence with the ASMs. Therefore if our
conjecture is true then one can investigate the FPL model, and the ASMs,
using the methods elaborated in mathematical physics for integrable model
(see, for example, \cite{Baxter}).

The background information on the ASMs and their different combinatorial
forms can be found in the recent review by Propp \cite{Propp} and in
references therein. For more details on enumerative problems related
to the states of the FPL model see the article by Wieland \cite{Wieland}.

Following the review paper by Propp \cite{Propp} we define the
``generalized tic-tac-toe'' graph as the graph formed by $n$ horizontal
lines and $n$ vertical lines meeting $n^2$ intersections of degree $4$,
with $4n$ vertices of degree $1$ at the boundary. Then we number the
vertices of degree $1$. We start with the left top vertex and number
clockwise every other vertex. Now consider subgraphs of the  underlying
tic-tac-toe graph such that each of the $n^2$ internal vertices lies on
exactly two of the selected edges and each numbered external vertex lies on
a selected edge, while each unnumbered external vertex does not lie on a
selected edge (see, for example, figure \ref{dpfl}).
\unitlength 1bp
\begin{figure}[ht]
  \centering
  \begin{minipage}[t]{.425\linewidth}
    \centering
    \begin{picture}(98,98)
      \put(5,5){\includegraphics{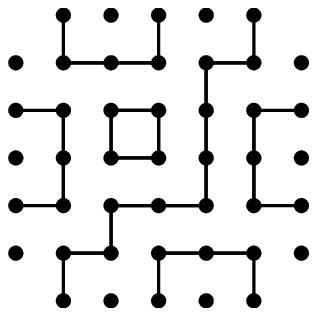}}
      \put(22,97){\makebox(0,0)[c]{\scriptsize 1}}
      \put(49,97){\makebox(0,0)[c]{\scriptsize 2}}
      \put(76,97){\makebox(0,0)[c]{\scriptsize 3}}
      \put(97,63){\makebox(0,0)[c]{\scriptsize 4}}
      \put(97,36){\makebox(0,0)[c]{\scriptsize 5}}
      \put(76,1){\makebox(0,0)[c]{\scriptsize 6}}
      \put(49,1){\makebox(0,0)[c]{\scriptsize 7}}
      \put(22,1){\makebox(0,0)[c]{\scriptsize 8}}
      \put(2,36){\makebox(0,0)[c]{\scriptsize 9}}
      \put(0,63){\makebox(0,0)[c]{\scriptsize 10}}
    \end{picture}
    \caption{One of the possible states of the FPL model for $n=5$}
    \label{dpfl}
  \end{minipage}
  \hspace{.05\linewidth}
  \begin{minipage}[t]{.425\linewidth}
    \centering
    \begin{picture}(98,98)
      \put(5,5){\includegraphics{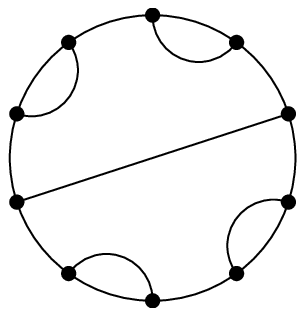}}
      \put(49,97){\makebox(0,0)[c]{\scriptsize 1}}
      \put(80,86){\makebox(0,0)[c]{\scriptsize 2}}
      \put(96,63){\makebox(0,0)[c]{\scriptsize 3}}
      \put(96,36){\makebox(0,0)[c]{\scriptsize 4}}
      \put(80,10){\makebox(0,0)[c]{\scriptsize 5}}
      \put(49,1){\makebox(0,0)[c]{\scriptsize 6}}
      \put(18,10){\makebox(0,0)[c]{\scriptsize 7}}
      \put(2,36){\makebox(0,0)[c]{\scriptsize 8}}
      \put(2,63){\makebox(0,0)[c]{\scriptsize 9}}
      \put(17,86){\makebox(0,0)[c]{\scriptsize 10}}
    \end{picture}
    \caption{The link-pattern corresponding to figure \ref{dpfl}}
    \label{lp}
  \end{minipage}
\end{figure}
These subgraphs are the states of the FPL model. The states the FPL model
are in bijective correspondence with the states of the square-ice model
with the domain-wall boundary conditions \cite{KIB}, and with the ASM
matrices. In particular, the number of such states is equal to the number
of ASMs, usually denoted by $A_n$. Each state of the FPL model define a
so-called link-pattern describing the pairings of the external vertices. We
depict such a pattern as a circle with $2n$ vertices placed on it and
connected pairwise inside the circle without intersection (see figure
\ref{lp}). Having in mind the relation to ASMs we denote the number of the
states of the FPL model corresponding to the link-pattern $\pi$ by
$A_n(\pi)$.

Consider as an example the case of $n=4$. In this case the FPL model has
42 states and 14 different link-patterns which are given in table
\ref{lps}.
\begin{table}[ht]
\caption{All possible link-patterns for $n=4$}
\label{lps}
\begin{center}
\begin{tabular}{lc}
\hline \hline
\multicolumn{1}{c}{\vrule height 1.1em width 0pt $\pi$} & $A_n(\pi)$ \\
\hline \\[-.8em]
\pairing{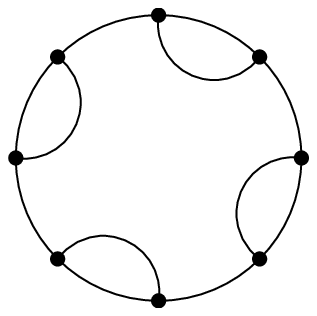}{.3} \hspace{1em} \pairing{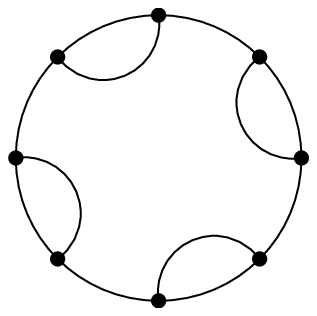}{.3} & \raisebox{12bp}{7} \\
\hline \\[-.8em]
\pairing{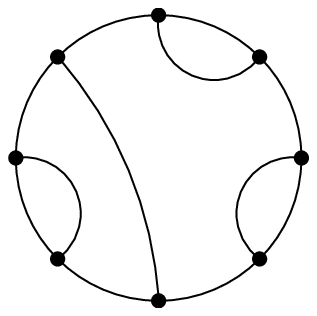}{.3} \hspace{1em} \pairing{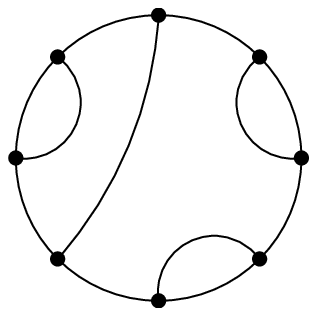}{.3} \hspace{1em} \pairing{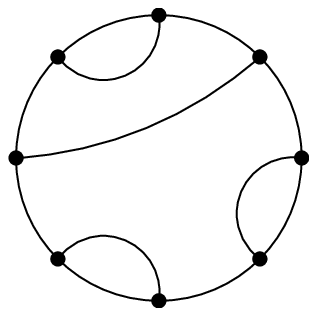}{.3}
\hspace{1em} \pairing{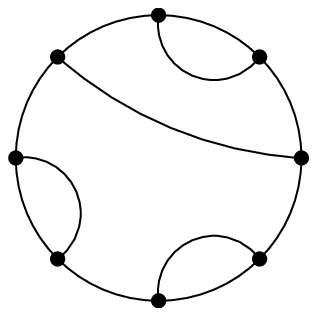}{.3} & \\[-.5em]
& 3 \\[-.5em]
\pairing{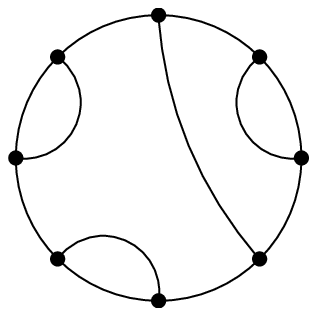}{.3} \hspace{1em} \pairing{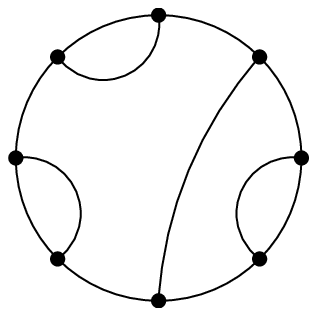}{.3} \hspace{1em} \pairing{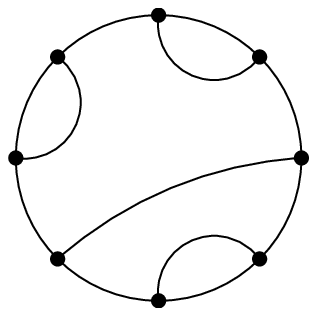}{.3}
\hspace{1em} \pairing{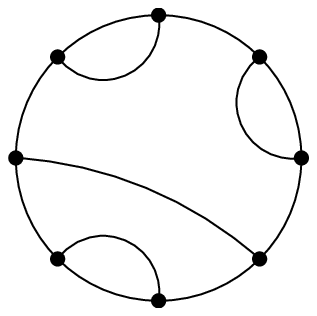}{.3} & \\
\hline \\[-.8em]
\pairing{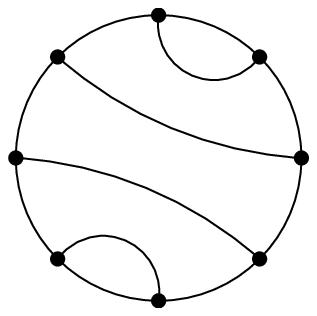}{.3} \hspace{1em} \pairing{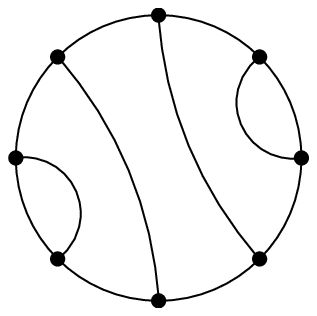}{.3} \hspace{1em}
\pairing{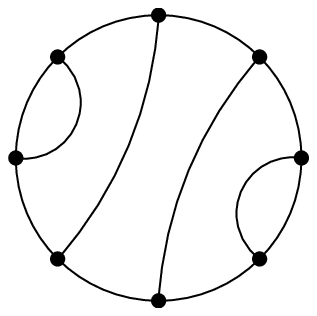}{.3} \hspace{1em} \pairing{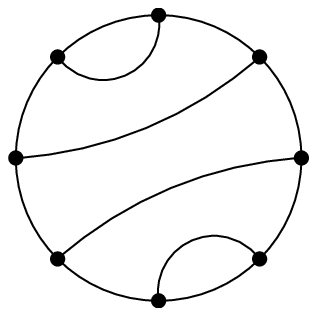}{.3} & \raisebox{12bp}{1} \\
\hline
\end{tabular}
\end{center}
\end{table}
The coincidence of the numbers $A_n(\pi)$ for different link-patterns is
explained by the conjecture of Bosley and Fidkowski, refined by
Cohn and Propp, and proved by Wieland \cite{Wieland}. In accordance
with this conjecture, if two link-patterns $\pi$ and $\pi'$ are connected
by a rotations or a reflection then $A_n(\pi) = A_n(\pi')$.

Let us now formulate a conjecture on the numbers $A_n(\pi)$. To this end
define $2n$ operations $h_i$, $i = 1,\dots 2n$, on link-patterns in the
following way. For a fixed $i$ consider a general link-pattern. Let this
link-pattern has the $i$-th vertex linked to the $j$-th one and the
$(i+1)$-th vertex linked to the $k$-th one. Two different case are
possible, the first is when $j = i+1$ and $k = i$. In this case the
operation $h_i$ left the link-pattern unchanged. In the second case $j \ne
i+1$ and $k \ne i$. In this case the operation $h_i$ replaces the two links
under consideration by the link of $i$-th vertex and the $(i+1)$-th vertex
and the link of the $j$-th vertex and the $k$-th vertex. To be
understandable by the widest audience of readers, we decided to formulate
our conjecture in three different ways.

The first formulation is in terms of the game theory. Let two players A and
B play a game of chance. Both persons randomly choose a state of the FPL
model. The player A checks a link-pattern of his state. The player B has to
choose randomly one of the  operations $h_i$ and to apply it to the
link-pattern of his state. Only the resulting link-pattern is of importance
for him. Players aim at obtaining some fixed link-pattern.
\begin{con}
Both players have equal chances to win.
\end{con}
Consider, for example, the case of $n=4$ with the wining link-pattern be
the very first one of table \ref{lps}. The chances of the player A are
$\mathsf P_A = 7/42 = 1/6$. Let us evaluate the chances of the player
B.
\begin{eqnarray*}
\mathsf P_{\mathrm B}  &=&
\mathsf P \left[\raisebox{-5bp}{\pairing{1}{.2}} \right] \mathsf P[\{h_1,
h_3, h_5, h_7\}] +
\mathsf P \left[\raisebox{-5bp}{\pairing{3}{.2}} \right] \mathsf P[\{h_5,
h_7\}] \\
&+& \mathsf P \left[\raisebox{-5bp}{\pairing{5}{.2}} \right]
\mathsf P[\{h_1, h_7\}] +
\mathsf P \left[\raisebox{-5bp}{\pairing{7}{.2}} \right] \mathsf P[\{h_1,
h_3\}] \\
&+& \mathsf P \left[\raisebox{-5bp}{\pairing{9}{.2}} \right] \mathsf
P [\{h_3, h_5\}] +
\mathsf P \left[\raisebox{-5bp}{\pairing{11}{.2}} \right] \mathsf P[\{h_3,
h_7\}] +
\mathsf P \left[\raisebox{-5bp}{\pairing{13}{.2}} \right] \mathsf P[\{h_1,
h_5\}].
\end{eqnarray*}
Using the data from table \ref{lps} one obtains
\[
\frac{7}{42}\times  \frac{4}{8}+\frac{3}{42}\times  \frac{2}{8}
+\frac{3}{42}\times  \frac{2}{8}+\frac{3}{42} \times \frac{2}{8}
+\frac{3}{42}\times  \frac{2}{8}+\frac{1}{42}\times  \frac{2}{8}
+\frac{1}{42}\times  \frac{2}{8}=\frac{1}{6}\,.
\]
Thus our conjecture in this case is true.

The second formulation is in terms of the combinatorics.
\begin{con} \label{con2}
For any $n = 1,2,\ldots$, one has
\[
\sum_{i=1}^{2n} \sum_{\pi' \in  \Pi_i(\pi)} A_n(\pi') = 2 n A_n(\pi),
\]
where  $\Pi_i(\pi)$ is the set of link-patterns $\pi'$ such that $h_i(\pi')
= \pi$.
\end{con}

The last formulation is in terms of the dense periodic O$(1)$ model. The
state space for this model is the set of all formal linear combinations of
the link-patterns. The operations $h_i$ extended to this space by linearity
become linear operators. The Hamiltonian of the model is the sum of the
operators $h_i$:
\begin{equation}
H = \sum_{i=1}^{2n} h_i. \label{ham}
\end{equation}
Consider the action of this operator on the vector
\begin{equation}
\Psi = \sum_{\pi} \pi A_n (\pi). \label{eig}
\end{equation}
We have
\[
H \Psi = \sum_{i=1}^{2n} \sum_\pi h_i(\pi) A_n(\pi) = \sum_\pi \pi
\sum_{i=1}^{2n}
\sum_{\pi' \in \Pi_i(\pi)} A_n(\pi') = 2 n \sum_\pi \pi A_n(\pi),
\]
where the last equality follows from conjecture \ref{con2}. Hence we can
write
\[
H \Psi = 2n \Psi.
\]
Thus, the vector $\Psi$ is an eigenvector of the Hamiltonian of the dense
periodic O$(1)$ loop model with an odd number of sites.

Let us consider our usual case of $n=4$. The matrix of the Hamiltonian $H$
and the component form of the vector $\Psi$ have in this case the form
\[
H = \left(\begin{array}{cccccccccccccc}
4 & 0 & 2 & 0 & 2 & 0 & 2 & 0 & 2 & 0 & 2 & 0 & 2 & 0 \\
0 & 4 & 0 & 2 & 0 & 2 & 0 & 2 & 0 & 2 & 0 & 2 & 0 & 2 \\
1 & 0 & 3 & 0 & 0 & 1 & 0 & 1 & 0 & 0 & 0 & 2 & 0 & 0 \\
0 & 1 & 0 & 3 & 0 & 0 & 1 & 0 & 1 & 0 & 0 & 0 & 2 & 0 \\
1 & 0 & 0 & 0 & 3 & 0 & 0 & 1 & 0 & 1 & 0 & 0 & 0 & 2 \\
0 & 1 & 1 & 0 & 0 & 3 & 0 & 0 & 1 & 0 & 2 & 0 & 0 & 0 \\
1 & 0 & 0 & 1 & 0 & 0 & 3 & 0 & 0 & 1 & 0 & 2 & 0 & 0 \\
0 & 1 & 1 & 0 & 1 & 0 & 0 & 3 & 0 & 0 & 0 & 0 & 2 & 0 \\
1 & 0 & 0 & 1 & 0 & 1 & 0 & 0 & 3 & 0 & 0 & 0 & 0 & 2 \\
0 & 1 & 0 & 0 & 1 & 0 & 1 & 0 & 0 & 3 & 2 & 0 & 0 & 0 \\
0 & 0 & 0 & 0 & 0 & 1 & 0 & 0 & 0 & 1 & 2 & 0 & 0 & 0 \\
0 & 0 & 1 & 0 & 0 & 0 & 1 & 0 & 0 & 0 & 0 & 2 & 0 & 0 \\
0 & 0 & 0 & 1 & 0 & 0 & 0 & 1 & 0 & 0 & 0 & 0 & 2 & 0 \\
0 & 0 & 0 & 0 & 1 & 0 & 0 & 0 & 1 & 0 & 0 & 0 & 0 & 2
\end{array}\right), \qquad
\Psi = \left( \begin{array}{c}
7 \\ 7 \\ 3 \\ 3 \\ 3 \\ 3 \\ 3 \\ 3 \\ 3 \\ 3 \\ 1 \\ 1 \\ 1 \\ 1
\end{array} \right).
\]
Here we take the basis vectors in the order given in table \ref{lps}.
Note that the sum of the matrix elements of $H$ belonging to any its fixed
column is equal to 8 which coincides with $2n$. This is a general property
of the Hamiltonian $H$ valid for any $n$ which implies that the spectral
radius of $H$ is equal to $2n$ (see, for example, \cite[section
9.3]{Lan69}). Therefore,
the vector $\Psi$ corresponds to the largest eigenvalue of the Hamiltonian
$H$. Actually this eigenvalue is nondegenerate that can be proved using the
positivity of the corresponding transfer matrix. Taking all this into
account we give the last formulation of our conjecture.

\begin{con}
The vector $(\ref{eig})$ is the unique eigenvector of the Hamiltonian
$(\ref{ham})$ corresponding to its largest eigenvalue $2n$.
\end{con}

In paper \cite{BGN} Batchelor, de Gier and Nienhuis formulated two
conjectures on the eigenvector under consideration\footnote{Actually
Batchelor, de Gier and Nienhuis take the Hamiltonian with the minus sign.
In this case the eigenvector under consideration is the ground state of the
model.}. They conjectured that the sum of its components is equal to $A_n$
and its largest component is equal to $A_{n-1}$. Our conjecture reveals the
combinatorial nature of all components of the eigenvector with the largest
eigenvalue.

Note that the Hamiltonian $H$ has the rotational and reflection invariance.
It can be easily shown that the eigenvector with the maximal eigenvalue
must have rotational and reflection invariance. Hence, our conjecture
generalizes partially the conjecture proved by Wieland in paper
\cite{Wieland} and mentioned above.

In conclusion note that we have verified our conjecture up to $n=7$.

{\it Acknowledgments} The authors would like to thank Bernard
Nienhuis for his clear explanation of the Hamiltonian for the dense O$(1)$
loop model. We have to acknowledge the beautiful survey of James Propp
which inspired us to investigate the FPL model. We are also thankful to him
for his interest to our work and discussions.

The work was supported in part by the Russian Foundation for Basic Research
under grant \#01--01--00201.

\end{document}